\documentclass[letterpaper,twoside,10pt]{article}
\usepackage{amsmath,amsthm,amssymb, enumerate}
             \newcommand{\marginalnote}[1]{}
\swapnumbers
\theoremstyle{plain}
\newtheorem{Thm}{Theorem}[section]
\newtheorem{Prop}[Thm]{Proposition}
\newtheorem{Lem}[Thm]{Lemma}
\newtheorem{Cor}[Thm]{Corollary}

\newtheorem*{Thm*}{Theorem}
\theoremstyle{remark}

\theoremstyle{definition}
\newtheorem{Rem}[Thm]{Remark}
\newtheorem{Rems}[Thm]{Remarks}

\newtheorem{Not}[Thm]{Notation}

\newtheorem{Defs}[Thm]{Definitions}
\newtheorem{Hyps}[Thm]{Hypotheses}
\numberwithin{equation}{Thm}
\newcommand{\abs}[1]{\left\lvert#1\right\rvert} 
\newcommand{\norm}[1]{\left\lVert#1\right\rVert} 
\newcommand{\gen}[1]{\left\langle#1\right\rangle}
\newcommand{\normgen}[1]{\langle \hskip -.07cm \lvert \hskip .05cm  #1 \hskip .05cm  \rvert \hskip -.07cm \rangle}
\def\d1{\discretionary{-}{}{-}}
\def\ors{sur\-face\d1plus\d1one\d1re\-la\-tion }
\def\Ors{Sur\-face\d1plus\d1one\d1re\-la\-tion }
\def\ore{one\d1re\-la\-tor }
\def\Ore{One\d1re\-la\-tor }

\DeclareMathOperator{\Hop}{H}
\DeclareMathOperator{\tr}{tr}

\DeclareMathOperator{\im}{im}
\DeclareMathOperator{\Tor}{Tor}
\DeclareMathOperator{\hd}{hd}
\DeclareMathOperator{\cd}{cd}
\DeclareMathOperator{\FL}{FL}
\DeclareMathOperator{\VFL}{VFL}
\DeclareMathOperator{\Z}{Z}

\DeclareMathOperator{\uE}{\underline{E}}
\begin{document}

\pagestyle{myheadings}
\markboth{$L^2$-Betti numbers of \ore  groups}{W. Dicks and P.\ A.\ Linnell}
\title{$L^2$-Betti numbers of \ore groups}

\author{Warren Dicks and  Peter A.\ Linnell}

\date{\today}

\maketitle

\begin{abstract}
We determine the $L^2$-Betti numbers of all \ore  groups and
all \ors groups.
We also obtain some information about the $L^2$-co\-homo\-logy of left\d1order\-able groups,
and deduce the non-$L^2$ result that, 
in any left\d1order\-able group of homological dimension one, 
all  two\d1gen\-er\-a\-tor subgroups are free.

\medskip

{\footnotesize

}
\end{abstract}

\maketitle

\section{Notation and background} \label{Sbackground}
Let $G$ be a (discrete) group, fixed throughout the article.

We use $\mathbb{R} \cup \{-\infty,\infty\}$
with the usual conventions; for example,
$\frac{1}{\infty} = 0$, and $3-\infty = -\infty$.
Let $\mathbb{N}$ denote the set of finite cardinals, $\{0,1,2,\ldots\}$.
We call $\mathbb{N} \cup \{\infty\}$ the 
set of {\it vague cardinals}, and, for each set $X$, we define its
 {\it vague cardinal} $\abs{X} \in \mathbb{N} \cup \{\infty\}$  
to be the cardinal of $X$ if $X$ is finite, and to be $\infty$ if $X$ is infinite.

Mappings of right modules will be
written on the left of their arguments, 
and mappings of left modules will be written on the right of their arguments.

Let $\mathbb{C}[[G]]$ denote the set of all functions from $G$ to $\mathbb{C}$
expressed as formal sums, that is, a function $a \colon G \to \mathbb{C}$, $g \mapsto a(g)$, 
will be written as $\sum_{g \in G} a(g) g$. Then $\mathbb{C}[[G]]$ 
has a natural $\mathbb{C}G$-bimodule structure, and contains
a copy of $\mathbb{C}G$ as $\mathbb{C}G$-sub-bimodule.
For each $a \in \mathbb{C}[[G]]$, 
we define $\norm{a} := (\sum_{g \in G} \abs{a(g)}^2)^{1/2} \in [0,\infty]$,
and $\tr(a) := a(1) \in \mathbb{C}$.

Define 
$$l^2(G) := \{a \in \mathbb{C}[[G]] :  \norm{a} < \infty\}.$$
We view $\mathbb{C} \subseteq \mathbb{C}G \subseteq l^2(G) \subseteq \mathbb{C}[[G]].$
There is a well-defined 
{\it external multiplication} map
$$l^2(G) \times l^2(G) \to \mathbb{C}[[G]], \quad 
(a ,b ) 
\mapsto a \cdot b,$$
where, for each $g \in G$,
$(a\cdot b)(g) := \sum_{h \in G} a(h) b(h^{-1}g)$;
this sum converges in~$\mathbb{C}$, and, moreover, 
 $\abs{(a\cdot b)(g)} \le \norm{a} \norm{b}$, by the Cauchy-Schwarz inequality.
The external multiplication extends the multiplication of $\mathbb{C}G$.

The {\it group von Neumann algebra of $G$},  denoted $\mathcal{N}(G)$, 
is the ring of bounded $\mathbb{C}G$-endo\-mor\-phisms of the
right $\mathbb{C}G$-module $l^2(G)$; see~\cite[\S1.1]{Lueck02}.  Thus
$l^2(G)$ is an $\mathcal{N}(G)$-$\mathbb{C}G$-bi\-mod\-ule.  We 
view $\mathcal{N}(G)$ as a subset of $l^2(G)$ by the map
$\alpha \mapsto \alpha(1)$, where $1$ denotes the identity element of 
$\mathbb{C}G \subseteq l^2(G)$.  
It can be shown that 
$$\mathcal{N}(G) = 
\{a \in l^2(G) \mid a \cdot l^2(G) \subseteq l^2(G)\},$$
and that the action of $\mathcal{N}(G)$ on $l^2(G)$ is given by 
the external multiplication.
Notice that $\mathcal{N}(G)$ contains  $\mathbb{C}G$ as a subring
and also that we have an induced `trace 
map' $\tr \colon \mathcal{N}(G) \to \mathbb {C}$.
The elements of $\mathcal{N}(G)$ which are injective, 
as operators on $l^2(G)$, are precisely
the (two-sided) non-zerodivisors in $\mathcal{N}(G)$, 
and they form a left and right Ore subset of 
$\mathcal{N}(G)$; see~\cite[Theorem~8.22(1)]{Lueck02}.

Let $\mathcal{U}(G)$ denote the 
{\it ring of unbounded operators affiliated to $\mathcal{N}(G)$}; 
see~\cite[\S8.1]{Lueck02}.   
It can be shown that $\mathcal{U}(G)$ is the left, and the right,
Ore localization of $\mathcal{N}(G)$ at the set of its non-zerodivisors.
For example, it is then clear that, 
\begin{equation}\label{Einvertible}
\text{if $x$ is an element of $G$ of infinite order, then $x-1$ is invertible in $\mathcal{U}(G)$.}
\end{equation}

Moreover, $\mathcal{U}(G)$ is a von Neumann regular ring in which one-sided 
inverses are two-sided inverses, and, hence,
one-sided zerodivisors are two-sided zerodivisors;
see~\cite[\S8.2]{Lueck02}.

There is a continuous, additive von Neumann dimension
that assigns to every left $\mathcal{U}(G)$-module $M$ 
a value
$\dim_{\mspace{2mu}\mathcal{U}(G)}M \in [0,\infty]$;
see Definition 8.28 and Theorem 8.29 of~\cite{Lueck02}.
For example, 
\begin{equation}\label{Eidempotent}
\text{if $e$ is an idempotent element of $\mathcal{N}(G)$, then
$\dim_{\mspace{2mu}\mathcal{U}(G)}\mathcal{U}(G)  e= \tr(e);$} 
\end{equation}
see Theorem 8.29 and \S\S6.1-2 of~\cite{Lueck02}. 

Consider any subring $Z$ of $\mathbb{C}$, and any resolution of $Z$ by projective,
or, more generally, flat, left $ZG$-modules
\begin{equation} \label{Eresolution}
 \cdots \longrightarrow P_2 \longrightarrow P_1
\longrightarrow P_0 \longrightarrow Z \longrightarrow 0,
\end{equation}
and let $\mathcal{P}$ denote the unaugmented complex
\begin{equation*} \label{Eprojectivecomplex}
 \cdots \longrightarrow P_2 \longrightarrow P_1
\longrightarrow P_0  \longrightarrow 0.
\end{equation*}
By Definition~6.50, Lemma~6.51 and Theorem 8.29 of~\cite{Lueck02}, 
we can define, for each $n \in \mathbb {N}$, 
the {\it $n$th $L^2$-Betti number of $G$} as
$$b_n^{(2)}(G) := \dim_{\mspace{2mu}\mathcal{U}(G)} 
\Hop_n(\mathcal{U}(G) \otimes_{ZG} \mathcal{P}),$$ 
where  $\mathcal{U}(G)$ is to be 
viewed as a $\mathcal{U}(G)$-$ZG$-bi\-module.
Of course, $$\Hop_n(\mathcal{U}(G) \otimes_{ZG} \mathcal{P})
=\Tor_n^{ZG}(\mathcal{U}(G), Z) \simeq  
\Tor_n^{\mathbb{Z}G}(\mathcal{U}(G), \mathbb{Z}) =  \Hop_n(G; \mathcal{U}(G)),$$
 where, for the purposes of this article,  
it will be convenient to understand that $\Hop_n(G; -)$
applies to {\it right} $G$-modules.
Thus the $L^2$-Betti numbers do not depend on the choice of $Z$, nor
on the choice of $\mathcal{P}$. 

\begin{Rem}\label{Rb0}
If $G$ contains an element of infinite
order, then~\eqref{Einvertible} implies that $ \mathcal{U}(G) \otimes_{ZG} Z= 0$,
and $ \mathcal{U}(G) \otimes_{ZG} P_1 \longrightarrow \mathcal{U}(G) \otimes_{ZG} P_0 
\longrightarrow 0$ is exact, and $\Hop_0(G;\mathcal{U}(G)) =  0$, and $b_0^{(2)}(G) = 0$.
\hfill\qed 
\end{Rem}

\begin{Rems}
In general, there is little relation between the $n$th
$L^2$-Betti number, $b_n^{(2)}(G) 
= \dim_{\mspace{2mu}\mathcal{U}(G)}\Hop_n(G; \mathcal{U}(G)) \in [0,\infty]$,
and the $n$th (ordinary)  Betti number, 
$$b_n(G) :=\dim_{\mathbb {Q}} \Hop_n (G;\mathbb {Q}) \in [0,\infty].$$

We say that $G$ is {\it of type $\FL$} if, for $Z = \mathbb{Z}$, 
there exists a resolution~\eqref{Eresolution}  
such that all the $P_n$ are finitely generated free left $\mathbb {Z}G$-modules
and all but finitely many of the $P_n$ are $0$.  

If $G$ is of type FL, then it is easy to see that 
the {\it $L^2$-Euler characteristic} 
$$\chi^{(2)}(G):= \sum\limits_{n \ge 0} (-1)^n b_n^{(2)}(G)$$ is equal to the
(ordinary)
{\it  Euler characteristic} $$\chi(G):=\sum\limits_{n \ge 0} (-1)^n b_n(G).$$

We say that $G$ is {\it of type $\VFL$} if $G$ has a subgroup $H$ of finite index
such that $H$ is of type $\FL$.  In this event, the (ordinary) {\it 
 Euler characteristic} of $G$ is
defined as $\chi(G) := \frac{1}{[G:H]}\chi(H)$; this is sometimes called the virtual
 Euler characteristic.  Here again, $\chi^{(2)}(G) = \chi(G)$; 
see~\cite[Remark~6.81]{Lueck02}. 
\hfill\qed
\end{Rems}

\section{Summary of results} \label{Ssummary}

In outline, the article has the following structure.
More detailed definitions can be found in the appropriate sections.

In Section~\ref{SU(G)}, we prove a useful technical result about
$\mathcal{U}(G)$ for special types of groups.

In Section~\ref{Sonerelator}, we calculate the $L^2$-Betti numbers of  
\ore  groups.  Let us describe the results.

For any element $x$ of a group $G$, we define the
{\it exponent} of $x$ in~$G$, denoted $\exp_G(x)$, as 
the supremum in $\mathbb{Z} \cup \{\infty\}$ 
of the set of those integers $m$ such that 
$x$ equals the $m$th power of some element of $G$. 
Then $\exp_G(x)$ is a nonzero vague cardinal.
We write $G/\normgen{x}$ to denote the quotient group of
$G$ modulo the normal subgroup of $G$ generated by $x$.

Suppose that $G$ has a one-relator
presentation $\gen{X\mid r}$. Thus $r$ is an element of the 
free group~$F$ on~$X$, and $G = F/\normgen{r}$. 

 Set $d := \abs{X} \in [0,\infty]$,
$m := \exp_F(r) \in [1,\infty]$, and
$\chi:= 1 - d + \frac{1}{m} \in [-\infty,1]$. 

It is known that
if $d < \infty$ then $G$ is of type VFL and
$\chi(G) = \chi$.  If $d = \infty$, then $G$ is not finitely generated and
$\chi = -\infty$; here we {\it define} $\chi(G) = -\infty$, which is non-standard, 
but it is reasonable.
 
In general,  $\max\{\chi(G),0\} = \frac{1}{\abs{G}}$.

In Theorem~\ref{Tonerelator}, we will show that,
\begin{equation}\label{Ebetti}
\text{for $n \in \mathbb{N}$, \quad}
b_n^{(2)}(G) = \begin{cases}
\max\{\chi(G),0\}  &\text{if $n = 0$,}\\
\max\{-\chi(G),0\} &\text{if $n = 1$,}\\
0 &\text{if $n \ge 2 $. }
\end{cases}
\end{equation}
L\"uck \cite[Example 7.19]{Lueck02} gave some results and conjectures 
concerning the $L^2$-Betti numbers of torsion-free 
\ore  groups, and~\eqref{Ebetti} shows
that the conjectured statements are true.

In Section~\ref{Sors}, we calculate the $L^2$-Betti numbers of  
an arbitrary \ors group  $G = \pi_1(\Sigma)/\normgen{\alpha}$.
Here $\Sigma$ is a connected orientable surface, and $\alpha$ is an
element of the fundamental group, $\pi_1(\Sigma)$.
The \ors groups were introduced and studied by
Hempel \cite{Hempel90}, and further investigated by Howie~\cite{Howie05};
these authors called the groups 
`one\d1re\-la\-tor surface groups',
but we are reluctant to adopt this terminology.

If $\Sigma$ is not closed, then $\pi_1(\Sigma)$ is a 
countable free group, see~\cite{Richards63}, and $G$ is a countable \ore  group.
In light of Theorem~\ref{Tonerelator}, 
we may assume that $\Sigma$ is a closed surface.  

Let $g$ denote the
genus of the closed surface $\Sigma$, and let $m = \exp_{\pi_1(\Sigma)}(\alpha)$.  
It is not difficult to deduce from known results that
$G$ is of type VFL and  $$\chi(G) = \begin{cases}
1 &\text{if $g = 0$,}\\
0 &\text{if $g = 1$,}\\
2 - 2g + \frac{1}{m}&\text{if $g \ge 2$.} 
\end{cases}$$
Then $\chi(G) \in (-\infty, 1]$ and
$\max\{\chi(G),0\} = \frac{1}{\abs{G}}$. In Section~\ref{Sors},
we will show that~\eqref{Ebetti} is also valid for \ors groups.

For any group $G$, $b_0^{(2)}(G)  = \frac{1}{\abs{G}}$; 
see~\cite[Theorem~6.54(8)(b)]{Lueck02}.
It is obvious that if $G$ is finite then $b_n^{(2)}(G) = 0$ for all $n \ge 1$.  
Thus, in essence, the foregoing results assert that
if $G$ is an infinite \ore  group, or an infinite \ors group, then
$$b_n^{(2)}(G) = \begin{cases}
-\chi(G) &\text{if $n = 1$,}\\
0 &\text{if $n \ne 1 $,}
\end{cases}$$
and we emphasize that, in this case, 
we understand that $\chi(G) = -\infty$ if $G$ is not finitely generated.

In Section~\ref{Sorderable}, we consider a variety of situations where
$Z$ is a nonzero ring and there exists some positive integer $n$ such that
 $P_n = ZG^2$ in a projective $ZG$-resolution~\eqref{Eresolution} of $\null_{_{ZG}}Z$.
For example, this happens for two-generator groups and for two-relator groups.

Thus, in Corollary~\ref{CThompson},  
we recover  L\"uck's result~\cite[Theorem 7.10]{Lueck02} that   
all the $L^2$-Betti numbers of Thompson's group $F$ vanish; see~\cite{CFP96} 
for a detailed exposition of the definition and main properties of $F$.   

\begin{Defs}
Recall that  $G$ is {\it left orderable} if there exists a
total order $\le$ of $G$ which is left $G$-invariant, 
that is, whenever $g,x,y \in G$ and $x \le y$,  
then $gx \le gy$.  One then says that $\le$ is a {\it left order} of $G$.
The reverse order is also a left order.
Since every group is isomorphic to its opposite through the inversion map,
we see that `left-orderable' is a short form for  `one-sided-orderable'.

A group is said to be {\it locally indicable} if every
finitely generated subgroup is either trivial or has an infinite cyclic
quotient. Burns and Hale~\cite{BurnsHale72} showed that every locally 
indicable  group is left orderable.  This often provides a
convenient way to prove that a given group is left orderable.

Recall that the {\it cohomological dimension of $G$ with respect to a ring $Z$}, 
denoted $\cd_{_{\scriptstyle Z}}G$, is the least $n \in \mathbb{N}$ such that 
$P_{n+1} = 0$ in some projective $ZG$-resolu\-tion~\eqref{Eresolution} 
of $_{_{\scriptstyle ZG}}Z$.  The {\it cohomological dimension} of $G$, denoted $\cd G$,
is $\cd_{_{\scriptstyle \mathbb{Z}}}G$.  A classic result of
 Stallings and Swan says that the groups of cohomological dimension 
at most one are precisely the free groups.

Similarly, the {\it homological dimension of $G$ with respect to a ring $Z$}, 
denoted $\hd_{_{\scriptstyle Z}}G$, is the least $n \in \mathbb{N}$ such that 
$P_{n+1} = 0$ in some flat $ZG$-resolution~\eqref{Eresolution} of 
$_{_{\scriptstyle ZG}}Z$.  The {\it homological dimension} of $G$, denoted $\hd G$,
is $\hd_{_{\scriptstyle \mathbb{Z}}}G$.
\hfill\qed
\end{Defs}

We understand that Robert Bieri, in the 1970's, 
first raised the question  as to whether
the groups of homological dimension at 
most one are precisely the locally free groups.
Notice that a locally free group has homological dimension at 
most one, since the augmentation ideal of a locally free group is a directed union of
finitely generated free left submodules.  Recently, in~\cite{KLL04}, 
it was proved that if the homological dimension of $G$ is at most one
and $G$ satisfies the Atiyah conjecture 
(or, more generally, the group ring $\mathbb{Z}G$
embeds in a one-sided Noetherian ring), 
then $G$ is locally free.
In Corollary~\ref{Chd1}, we show that if $G$ is locally indicable,
or, more generally, left orderable, and  
the homological dimension of $G$ is at most one,
then every {\it two-generator} subgroup of $G$ is free.

Finally, in Proposition~\ref{Ptwotwo}, we calculate the first three
$L^2$-Betti numbers
of an arbitrary left-orderable two-relator group of cohomological dimension
at least three.

\begin{Not}\label{Nmatrices}
We will frequently consider maps between free modules over a ring~$U$,
and we will use the following format.  

Let $X$ and $Y$ be sets. 

By an {\it $X \times Y$ row-finite matrix} over $U$
we mean a function $(u_{x,y}) \colon X \times Y \to U$, $(x,y) \mapsto u_{x,y}$ such 
that, for each $x \in X$, $\{y \in Y \mid u_{x,y} \ne 0\}$ is finite.

We write $\oplus_X U$ to denote the direct sum of copies of $U$ indexed by $X$.
If $n \in \mathbb{N}$, and $X = \{1,\ldots,n\}$, we identify $X = n$ and 
also write $\oplus_n U$ as $U^n$.
An element of $\oplus_X U$ will be viewed as a $1 \times X $ row-finite matrix
$(u_{1x})$ over $U$.  Then $\oplus_X U$ is a left $U$-module in a natural way.

A map  $\oplus_X U \to \oplus_Y U$ of left $U$-modules will be thought of as 
right multiplication by a row-finite $X \times Y$ matrix $(u_{x,y})$ in
a natural way, and we will write\linebreak $\oplus_X U \xrightarrow{(u_{x,y})} \oplus_Y U$.
\hfill\qed
\end{Not}

\section{Preliminary results about $\mathcal{U}(G)$} \label{SU(G)}

For $a = \sum_{g\in G} a(g)g \in \mathbb {C}[[G]]$, 
we let $a^* = \sum_{g\in G} \overline{a(g^{-1})}g$ 
where $\overline{z}$ indicates the complex
conjugate of $z$.  This involution restricts to
$\mathbb{C}(G)$ and $\mathcal{N}(G)$, and  
extends in a unique way to  $\mathcal{U}(G)$. 
Furthermore, if $a,
b \in \mathcal{N}(G)$, then $(a b)^* =
b^*a^*$ and $a^* a= 0$ if and only if $a = 0$.

In Sections 4 and 5, we shall see that the narrow hypotheses of the 
following result hold whenever $G$ is a one-relator group or a \ors group.  

\begin{Thm}\label{Text}  Suppose that $G$ has a normal subgroup $H$ 
such that $H$ is the semidirect product  $F \rtimes C$ 
of a free subgroup $F$ by a finite subgroup $C$,
and that $G/H$ is locally indicable, or, more generally, left orderable.

Let $m = \abs{C}$, and let $e = \frac{1}{m} \sum\limits_{c\in C} c  \in \mathbb{C}G$. 

Then the following hold.

\begin{enumerate}[\normalfont(i)]
\vskip-0.6cm \null
\item\label{Ione}  Each torsion subgroup of $G$ embeds in $C$. 
\vskip-0.6cm \null
\item\label{Itwo} Each nonzero element of $e\mathbb{C}Ge$ 
 is invertible in $e\mspace{2mu}\mathcal{U}(G)e$.
\vskip-0.6cm \null
\item\label{Ithree} For all $x  \in \mathcal{U}(G)e$ and $y \in e\mspace{2mu}\mathbb{C}G$, 
 if $xy = 0$ then $x = 0$ or $y = 0$.

\end{enumerate}

\end{Thm}

\begin{proof} \eqref{Ione}  Each torsion subgroup of $G$ lies in $H$ 
and has trivial intersection with $F$, and therefore embeds in $C$.

\eqref{Itwo}  Notice that $e$ is a projection, that is, $e$ is idempotent and 
$e^\ast= e$.  Clearly, $\tr(e) = \frac{1}{m}$.
Also, $e\mspace{2mu} \mathcal{U}(G)e$ is a ring and $e\mspace{2mu}\mathbb{C}Ge$ is a subring of 
$e\mspace{2mu} \mathcal{U}(G)e$.
Moreover, in $e \mspace{2mu}\mathcal{U}(G)e$, one-sided inverses  are two-sided inverses.

Let $a \in e\mathbb{C}Ge - \{0\}$.  We want to show that
$a$ is left invertible in $e \mspace{2mu}\mathcal{U}(G)e$.

Let  $T$ be a transversal for the right (or left) $H$-action on $G$, and suppose 
that $T$ contains $1$.  Write $a = t_1a_1 + \cdots + t_n a_n $
where the $t_i$ are distinct elements of  $T$, and, for each $i$,
 $a_i \in \mathbb{C}(H)e - \{0\}$.

Let $\preceq$ be a left order for~$G/H$.
We may assume that $t_1H \prec \cdots \prec t_nH$.
To show that $a$ is left invertible in $e\mspace{2mu}\mathcal{U}(G)e$,
it suffices to show that $(ea_1^\ast t_1^{-1}e)a$ 
is left invertible in $e\mspace{2mu}\mathcal{U}(G)e$.  On replacing $a$ with
 $(ea_1^*t_1^{-1}e)a = a_1^*t_1^{-1}a$,
we see that we may assume that $t_1 = 1$ and  $a_1  \in e\mathbb{C}He - \{0\}$.

By~\eqref{Ione}, $m$ is the least common multiple of the orders of the finite
subgroups of~$H$.  Now the strong Atiyah conjecture holds for $H$; 
see~\cite{Linnell93} or~\cite[Chapter 10]{Lueck02}.
Hence $\dim_{\mspace{2mu}\mathcal{U}(H)}\mathcal{U}(H)a_1 \ge \frac{1}{m} = \tr(e)$.
Of course, $\mathcal{U}(H)a_1 \subseteq \mathcal{U}(H)e$, and thus
$\dim_{\mspace{2mu}\mathcal{U}(H)}\mathcal{U}(H)a_1 
\le \dim_{\mspace{2mu}\mathcal{U}(H)}\mathcal{U}(H)e = \tr(e)$.
Hence $\dim_{\mspace{2mu}\mathcal{U}(H)}\mathcal{U}(H)a_1  = \tr(e)$.

Also, $\mathcal{U}(H)(a_1 + 1-e) = \mathcal{U}(H)a_1 \oplus  \mathcal{U}(H)(1-e)$.
Hence 
\begin{align*}
\dim_{\mspace{2mu}\mathcal{U}(H)}\mathcal{U}(H)(a_1 + 1-e)
&= \dim_{\mspace{2mu}\mathcal{U}(H)} \mathcal{U}(H)a_1 + \dim_{\mspace{2mu}\mathcal{U}(H)} \mathcal{U}(H)(1-e)\\
&= \tr(e) + \tr(1-e) = 1.
\end{align*}
This implies that $a_1 + 1-e$ is invertible in $\mathcal{U}(H)$.
The $\ast$-dual of~\cite[Theorem 4]{Linnell92} now implies that
$a + 1-e = 1(a_1 + 1-e) + t_2a_2 + \cdots + t_n a_n$ 
is invertible in $\mathcal{U}(G)$.  It is then straightforward to show that   
$a$ is invertible in $e\mspace{2mu} \mathcal{U}(G)e$.

\eqref{Ithree}  Suppose that $y \ne 0$.  Then $x^*xyy^* = 0$, 
$yy^* \in  e\mathbb{C}Ge - \{0\}$ and $x^*x \in e\mspace{2mu}\mathcal{U}(G)e$.
By~\eqref{Itwo}, $yy^*$ is invertible in $e\mspace{2mu}\mathcal{U}(G)e$.
Hence $x^*x = 0$ and $x = 0$.
\end{proof}

\begin{Rem} The above proof shows that the conclusions
of Theorem~\ref{Text}\eqref{Itwo} and~\eqref{Ithree} hold under the following
hypotheses: $H$ is a normal subgroup of $G$; 
$G/H$ is left orderable;  the
 strong Atiyah conjecture holds for $H$; and,  
$e$ is a nonzero projection in $\mathbb{C}H$ such that 
$\frac{1}{\tr(e)}$
is the least common multiple of the orders of the finite subgroups of~$H$.
\hfill\qed
\end{Rem}

The degenerate 
case of Theorem~\ref{Text}\eqref{Itwo} where $H = F = C =1$ 
follows directly
from~\cite[Theorem~2]{Linnell92}. 

\begin{Thm}\label{Torder}
If $G$ is locally indicable, or, more generally,  
left orderable, then
every nonzero element of $\mathbb{C}G$ is invertible in $\mathcal{U}(G)$.\hfill\qed
\end{Thm}

\section{\Ore groups} \label{Sonerelator}

We shall now calculate the $L^2$-Betti numbers of \ore  groups.

\begin{Not}\label{Nonerelator}  
Suppose that $G$ is a \ore  group, and let
$\gen{ X \mid r }$ be a  \ore presentation of $G$.

Here $r$ is an element of the free group $F$ on $X$ and $G = F/\normgen{r}$.

Let $m =\exp_F(r)$ and let $d = \abs{X}$.  
These are vague cardinals.
Here $m \ne 0$; moreover, $m = \infty$ if and only if $r=1$,
in which case $G = F$.

If $ m < \infty$, then $r = q^m$ for some $q \in F$.
Let $c$ denote the image of $q$ in $G$, and let $C = \gen{c} \le G$.
Then $C$ has order $m$.
Let $e = \frac{1}{m} \sum_{x \in C} x \in \mathbb{C}G$.

If $m = \infty$, we define $e = 0 \in \mathbb{C}G$.

In any event $e$ is a projection and $\tr(e) = \frac{1}{m}$.

There is an exact sequence of left $\mathbb{Z}G$-modules
\begin{align*}
0 \longrightarrow \oplus_X \mathbb{Z}G \longrightarrow \mathbb{Z}G \longrightarrow \mathbb{Z}  
\longrightarrow 0&&&\text{if $m = \infty$,}\\
0 \longrightarrow \mathbb{Z}[G/C] \longrightarrow \mathbb{Z}  
\longrightarrow 0&&&\text{if $d = 1$ and $m < \infty,$}\\
0  \longrightarrow \mathbb{Z}[G/C]
\longrightarrow \oplus_X \mathbb{Z}G \longrightarrow \mathbb{Z}G \longrightarrow \mathbb{Z}  
\longrightarrow 0&&&\text{if $d \ge 2$ and $m < \infty$.}
\end{align*}
see~\cite{Chiswell03}, specifically,  Lemma~6.21 and ($*$) on
 p.~167 in the proof of 
Theorem~6.22. 
In all cases, there is then an exact sequence of left $\mathbb{C}G$-modules
\begin{equation}\label{EChiswellC}
0 \longrightarrow  \mathbb{C}Ge \xrightarrow{(a_{1,x})} \oplus_{X} \mathbb{C}G 
\xrightarrow{ (b_{x,1})} \mathbb{C}G \longrightarrow \mathbb{C} \longrightarrow 0;
\end{equation}
for each $x \in X$,  $b_{x,1}$ 
is the image of $x-1$ in $\mathbb{C}G$, and $a_{1,x}$ is the left Fox derivative
$\frac{\partial r}{\partial x} = (me)\frac{\partial q}{\partial x}\in e\mathbb{C}G$.

If $d < \infty$, then $G$ is of type VFL and
\begin{equation}\label{Eeuler}
\chi(G) = 1 - d + \frac{1}{m} \in (-\infty,1];
\end{equation}
see Theorem~6.22  and Corollary~6.15 of~\cite{Chiswell03}, for the cases
where $m < \infty$ and $m = \infty$, respectively.

In the case where $d = \infty$, that is, 
$G$ is a non-finitely-generated \ore  group,
we {\it define} $\chi(G):= -\infty$.  This 
is non-standard, but it extends~\eqref{Eeuler}.

It is  easy to verify that
$\frac{1}{\abs{G}} = \max\{\chi(G),0\}$.   
In fact, by abelianizing $G$, we see that 
$G$ is finite if and only if
either $d = 1$ and $m < \infty$, or $d = 0$ (and hence $m = \infty$).  
\hfill\qed
\end{Not}

We shall now prove the following.

\begin{Thm} \label{Tonerelator}
If $G$ is a \ore  group, then, for $n \in \mathbb{N}$,
\begin{equation}\label{Echar}
b_n^{(2)}(G) = \begin{cases}
\max\{\chi(G),0\} \,\,(= \frac{1}{\abs{G}}) &\text{if $n = 0$,}\\
\max\{-\chi(G), 0\} &\text{if $n = 1$,}\\
0 &\text{if $n \ge 2 $. }
\end{cases}
\end{equation}
\end{Thm}

\begin{proof} Suppose that Notation~\ref{Nonerelator} holds.

Unaugmenting~\eqref{EChiswellC} and 
applying $\mathcal{U}(G) \otimes_{\mathbb{C}G} -$ gives
\begin{equation}\label{EChiswellU}
0 \longrightarrow  \mathcal{U}(G)e \xrightarrow{(a_{1,x})} \oplus_{X}\mathcal{U}(G)
\xrightarrow{ (b_{x,1})} \mathcal{U}(G)\longrightarrow  0;
\end{equation}
the homology of~\eqref{EChiswellU} is then $\Hop_\ast(G;\mathcal{U}(G))$.

We claim that
\begin{equation}\label{Eonerelator}
\text{if $y \in \mathcal{U}(G)e -\{0\}$ and
$a \in e\mathbb{C}G -\{0\}$, then $ya \ne 0$.}
\end{equation}

This is vacuous if $m = \infty$.

If $m < \infty$, let $H$ denote the normal subgroup of $G$ generated by $c$.
Then $G/H = \gen{X \mid q}$ is a torsion-free \ore  group.
Hence $G/H$ is locally indicable by 
\cite[Theorem 3]{Brodskii80}, \cite[Theorem 4.2]{Howie82} 
or \cite[Corollary 3.2]{Howie00}.
Also $H$ is the free product of
certain $G$-conjugates of $C$, by~\cite[Theorem 1]{FKS72}. 
 By mapping each of
these conjugates of $C$ isomorphically to $C$, we obtain an
epimorphism $H\twoheadrightarrow C$.
Applying \cite[Proposition I.4.6]{DicksDunwoody89} to this
epimorphism, we see that its kernel $F$ is free.  
Clearly, $H = F \rtimes C$.
Now~\eqref{Eonerelator} holds by Theorem~\ref{Text}\eqref{Ithree}.

Since $(a_{1,x})$ is injective in~\eqref{EChiswellC}, either $e = 0$ or
there is some $x_0 \in X$ such that 
$a_{1,x_0} \ne 0$.  
It follows from~\eqref{Eonerelator} that $(a_{1,x})$ is 
injective in~\eqref{EChiswellU}, and hence
$\Hop_2(G;\mathcal{U}(G)) = 0$.  
On taking $\mathcal{U}(G)$-dimensions, we find
$b_2^{(2)}(G) = 0$, and $\dim_{\mspace{2mu}\mathcal{U}(G)} \im ((a_{1,x})) = \frac{1}{m}$.

If either $d \ge 2$, or $d = 1$ and $m = \infty$, then,
by abelianizing, we see that  
there is some $x_1 \in X$ whose image in $G$ has infinite order.
By~\eqref{Einvertible}, we see that 
$(b_{x,1})$ is surjective in~\eqref{EChiswellU}, 
and hence $\Hop_0(G; \mathcal{U}(G))~=~0$.
On taking $\mathcal{U}(G)$-dimensions, we find that $b_0^{(2)}(G) = 0$,
$\dim_{\mspace{2mu}\mathcal{U}(G)} \im ((b_{x,1})) = 1$, and
$\dim_{\mspace{2mu}\mathcal{U}(G)} \ker ((b_{x,1})) = d-1$.
Now $$b_1^{(2)}(G) = \dim_{\mspace{2mu}\mathcal{U}(G)} \ker ((b_{x,1})) - 
\dim_{\mspace{2mu}\mathcal{U}(G)} \im ((a_{1,x})) = d-1-\frac{1}{m} = -\chi(G).$$  
Thus~\eqref{Echar} holds.

This leaves the cases where either $d = 0$ or $d = 1$ and $m < \infty$.
Here $G$ is finite cyclic, and again ~\eqref{Echar} holds.
\end{proof}

\section{\Ors groups} \label{Sors}

We next calculate the $L^2$-Betti numbers for an arbitrary \ors
group $G = \pi_1(\Sigma)/\normgen{\alpha}$, 
where $\Sigma$ is a connected orientable surface, possibly
with boundary and not necessarily compact, and $\normgen{\alpha}$ 
is the normal
closure of a single element $\alpha \in \pi_1(S)$.

By the results of the previous section, we may assume
that the implicit presentation of $G$ has more than one relator.
As explained in Section~\ref{Ssummary}, 
$\Sigma$ must a closed surface.  Let $g$ denote the genus of
$\Sigma$.  Then $g \in \mathbb{N}$ and  
\begin{equation*}\label{EpresPi}
\pi_1(\Sigma) =\gen{ \quad x_1,x_2,\ldots, x_{2g-1},x_{2g}  \quad \vert \quad
[x_1,x_2][x_3,x_4] \cdots [x_{2g-1}, x_{2g}] \quad },
\end{equation*}
where $[x,y]$ denotes $xyx^{-1}y^{-1}$.
Since this is a one-relator presentation, we have $\alpha \ne 1$.
In particular, $g$ is nonzero.  The 
non \ore cases are included in the following.

\begin{Thm}\label{Tors}  Let $\Sigma$ be a closed  orientable surface of genus at least one,
let\linebreak $S = \pi_1(\Sigma)$, let $\alpha$ be a nontrivial element of $S$, 
and let   $G = S/\normgen{\alpha}$. 

Let $g$ denote the genus of $\Sigma$, let $m = \exp_{S}(\alpha)$, 
and let $Q$ be a nonzero ring in which $\frac{1}{m}$ is defined, 
that is, if  $m < \infty$ then $mQ = Q$.
 Then the following hold.
\begin{enumerate}[\normalfont (i)]
\vskip-0.5cm \null
\item\label{Ibig1} $G$ is of type $\VFL$ and $\chi(G) = \min\{2-2g+\frac{1}{m},0\}=
\begin{cases}
0 &\text{if $g=1$,}\\
2 - 2g + \frac{1}{m} &\text{if $g \ge 2$.}
\end{cases}$
\vskip-0.3cm \null
\item\label{Ibig2} $\cd_QG = \min\{2,g\} = \begin{cases}
1 &\text{if $g=1$,}\\
2 &\text{if $g \ge 2$.}
\end{cases}$
\vskip-0.2cm \null
\item\label{Ibig3} For $n \in \mathbb{N}$,
$b_n^{(2)}(G) = -\delta_{n,1}\chi(G) = \begin{cases}
-\chi(G)  &\text{if $n = 1$,}\\
0 &\text{if $n \ne 1$. }
\end{cases} $
\end{enumerate}
\end{Thm}

\begin{proof}
We break the proof up into a series of lemmas and 
summaries of notation.

\begin{Not}
As in~\cite[Examples~I.3.5(v)]{DicksDunwoody89}, 
the expression $S_{1}\ast_{S_0} s$ will denote an HNN extension,
where it is understood that $S_1$ is a group, 
$S_0$ is a subgroup of $S_1$  and $s$ is  an injective group homomorphism
 $s \colon S_0 \to S_1$,  $a \mapsto a^s$.  The image of this homomorphism
 is denoted $S_0^s$.
\hfill\qed
\end{Not}

\begin{Lem}[Hempel]\label{LHempel}
If $g \ge 2$, then there exists an HNN-decomposition\linebreak 
$S = S_{1}\ast_{S_0} s$  such that $S_1$ is a free group,  
$\alpha$ lies in $S_1$,  and the normal subgroup of $S_1$ generated by
$\alpha$ intersects both $S_0$ and $S_0^s$ trivially.

Hence,   $G = S/\normgen{\alpha}$ has a matching HNN-decomposition 
$S/\normgen{\alpha} = S_{1}/\normgen{\alpha}\ast_{S_0} s$.
\end{Lem}

\begin{proof} This was implicit in the proof of \cite[Theorem 2.2]{Hempel90}, 
and was made explicit in \cite[Proposition 2.1]{Howie05}.
\end{proof}

\begin{Lem}[Hempel] If $m < \infty$, there exists $\beta \in S$ such that $\beta^m = \alpha$, 
and the image of $\beta$ in $G$ has order $m$.
\end{Lem}

\begin{proof} As this is obvious for $g = 1$, we may assume that $g \ge 2$.
Thus we have matching HNN-decompositions 
$S = S_{1}\ast_{S_0} s$ and
$G = S/\normgen{\alpha} = S_{1}/\normgen{\alpha}\ast_{S_0} s$,
as in Lemma~\ref{LHempel}.

Let $m' = \exp_{S_{1}}\alpha$.  Since $\alpha \ne 1$ and $S_1$ is free,
we see that $m' < \infty$.  Choose $\beta \in S_{1}$ 
such that $\beta^{m'} = \alpha$.  Let $c$ denote the 
image of $\beta$ in $G$, and let $C = \gen{c} \le G$.
Then  $C$ has order $m'$, and every torsion 
subgroup of $S_{1}/\normgen{\alpha}$ embeds in
$C$. From the HNN decomposition for $G$, we see that any finite 
subgroup of $G$ is conjugate to a subgroup of 
$S_{1}/\normgen{\alpha}$, and hence has order dividing $m'$. 

A similar argument shows that for any positive integer $i$,
$S/\normgen{\alpha^i}$ has a matching HNN decomposition, 
and therefore has a subgroup of order $m'i$ and a subgroup 
of order $i$.  It follows that if $\alpha = \gamma^j$ 
for some positive integer $j$ then $S/\normgen{\alpha}$ has a 
subgroup of order $j$, and hence $j$ divides $m'$.  
It now follows that $m = m'  < \infty$.
\end{proof}

\begin{Not}\label{NH}  Let $\beta$ denote an element of $S$ such that $\beta^m = \alpha$.

Let $c$ denote the image of $\beta$ in $G$.  Let $C = \gen{c}$, a cyclic 
subgroup of $G$ of order~$m$.  Let $e = \frac{1}{m}\sum_{x\in C} x$, an
idempotent element of $\mathbb{C}G$ with $\tr(e) = \frac{1}{m}$; we shall
also view $e$ as an idempotent element of $QG$. 

Let $H$ denote the normal subgroup of $G$ generated by $c$; thus,
$G/H~\simeq~S/\normgen{\beta}$.
\hfill\qed
\end{Not}

\begin{Lem}\label{L1} \begin{enumerate}[\normalfont(i)]
\vskip-0.6cm \null
\item\label{I3} $H$ has a free subgroup $F$ such that $H = F \rtimes C$.
\vskip-0.6cm \null
\item\label{I4} $G/H$ is locally indicable.  
\vskip-0.6cm \null
\item\label{I5} Every torsion subgroup of $G$ embeds in $C$.
\vskip-0.6cm \null
\item\label{I6} If $x \in \mathcal{U}(G)e -\{0\}$ and $y \in e\mathbb{C}G - \{0\}$, 
then $xy \ne 0$.
\end{enumerate}
\end{Lem}

\begin{proof} \eqref{I3}. As this is clear for $g=1$, we may assume that $g \ge 2$.

By Lemma~\ref{LHempel} with $\beta$ in place of $\alpha$,
there exists an HNN-decomposition\linebreak $S = S_{1}\ast_{S_0} s$ where
$S_1$ is a free group,  $\beta$ lies in $S_1$, and the normal subgroup 
of $S_1$ generated by $\beta$ intersects both $S_0$ and $S_0^s$ trivially.
Hence $\alpha$ lies in $S_1$, and the normal subgroup of $S_1$
generated by $\alpha$ intersects both $S_0$ and $S_0^s$ trivially.
It follows that we can make identifications
$$G = S/\normgen{\alpha} = S_{1}/\normgen{\alpha}\ast_{S_0} s \text{\quad and\quad}
G/H = S/\normgen{\beta} = S_{1}/\normgen{\beta} \ast_{S_0}s.$$
Thus we have matching HNN-decompositions for $S$, $G$ and $G/H$.

Let us apply  Bass-Serre theory, following, for 
example,~\cite[Chapter 1]{DicksDunwoody89}.
Consider the action of $H$ on the Bass-Serre tree for 
the above HNN-decomposition of $G$.  Then $H$ acts freely on the edges.
Let $H_0$ denote the normal subgroup of $S_{1}/\normgen{\alpha}$ 
generated by $c$.  Then $H_0$ is a
vertex stabilizer for the $H$-action, and the other
vertex stabilizers are $G$-conjugates of $H_0$.  
By Bass-Serre theory, or the Kurosh Subgroup Theorem, 
$H$ is the free product of a free group and various 
$G$-conjugates of $H_0$.

By~\cite[Theorem 1]{FKS72}, $H_0$ itself is a free product of certain
$S_1/\normgen{\alpha}$-con\-ju\-gates of~$C$.  

Thus $H$ is the free product a free group and 
various $G$-conjugates of $C$.    
If we map each of these $G$-conjugates of $C$
 isomorphically to $C$,  and map the free group to $1$,
we obtain an
epimorphism $H\twoheadrightarrow C$.
Applying \cite[Proposition I.4.6]{DicksDunwoody89} to this
epimorphism, we see that its kernel $F$ is free.  
Clearly, $H = F \rtimes C$. This proves~\eqref{I3}.

\eqref{I4}. Since $G/H = S/\normgen{\beta}$ and $\beta$ is not a proper
power in $S$, $G/H$ is locally indicable 
by~\cite[Theorem~2.2]{Hempel90}.  

\eqref{I5} and~\eqref{I6} hold by Theorem~\ref{Text}.
\end{proof}

Let us dispose of the case where $g=1$, which is well known and included
only for completeness.

\begin{Lem}\label{Lg1} If $g=1$, then the following hold. 
\begin{enumerate}[\normalfont(i)]
\vskip-0.6cm \null
\item $H = C$ and $G/C$ is infinite cyclic generated by $xC$ for some $x \in G$.
\vskip-0.6cm \null
\item $0 \longrightarrow \mathbb{Z}[G/C]  \xrightarrow{x-1}\mathbb{Z}[G/C] 
 \longrightarrow \mathbb{Z} \longrightarrow 0$
is an exact sequence of left $\mathbb{Z}G$-mod\-ules. 
\vskip-0.6cm \null
\item
$0 \longrightarrow QGe  \xrightarrow{x-1}QGe
 \longrightarrow Q \longrightarrow 0$ is an exact sequence of left $QG$-modules. 
\vskip-0.6cm \null
\item\label{Icd} $\gen{x}$ is an infinite cyclic subgroup of $G$ of finite index, 
$G$ is of type $\VFL$, $\chi(G)~=~0$ and $\cd_QG = 1$.
\vskip-0.6cm \null
\item  The homology of 
$0 \longrightarrow \mathcal{U}(G)e 
 \xrightarrow{x-1} \mathcal{U}(G)e  \longrightarrow 0$ is $\Hop_*(G;\mathcal{U}(G))$.
\vskip-0.6cm \null
\item\label{Ibs} For each $n \in \mathbb{N}$, $b_n^{(2)}(G) = 0$. \hfill\qed
\end{enumerate}
\end{Lem}

\begin{Rem} 
For $g = 1$, Lemma~\ref{Lg1}(ii) gives the augmented cellular chain
complex of a one-di\-men\-sion\-al $\uE(G)$ which resembles the real line.
\hfill\qed
\end{Rem}

\begin{Not}\label{not:der}  Henceforth we assume that $g \ge 2$.
 
Let $X = \{x_1,x_2,\ldots,x_{2g-1},x_{2g}\}$, let $F$ be the free group on $X$,
and let\linebreak $r_1 = [x_1,x_2]\cdots[x_{2g-1},x_{2g}] \in F.$
Then $S = \gen{X \mid r_1}$.  

Let  $q_2$ be any element of $F$ 
which maps to $\beta$ in $S$, and let $r_2 = q_2^m$.  
Then $G = \gen{X \mid r_1, r_2}$.

For $i \in \{1,2\}$, $j \in \{1,\ldots, 2g\}$, we set
$a_{i,j}:= \frac{\partial r_i}{\partial x_j} \in \mathbb{Z}G$,
the left Fox derivatives, and $b_{j,1}:= x_j -1 \in \mathbb{Z}G$. 

Notice that $me = \sum_{x \in C} x \in \mathbb{Z}G$ and
$a_{2,j} = \frac{\partial r_2}{\partial x_j}
 = (me) \frac{\partial q_2}{\partial x_j}$.
\hfill\qed
\end{Not}

\begin{Lem}[Howie]\label{L2} 
The sequence of left $\mathbb{Z}G$-modules
\begin{equation}\label{EL2}
 0 \longrightarrow
\mathbb{Z}G \oplus \mathbb{Z}[G/C] \xrightarrow{(a_{i,j})} 
\mathbb{Z}G^{2g} \xrightarrow{(b_{j,1})}
\mathbb{Z}G \longrightarrow \mathbb{Z} \longrightarrow 0
\end{equation}
is exact.
\end{Lem}

\begin{proof} Howie~\cite[Theorem~3.5]{Howie05} describes a
$K(G,1)$, and it is straightforward to give it a CW-structure as follows.

We take a $K(S,1)$ with  one zero-cell, $2g$ one-cells, and a 
two-cell which is a $2g$-gon,  and then the exact sequence of left
 $\mathbb{Z}S$-modules arising from the  augmented cellular 
chain complex of the 
universal cover of the $K(S,1)$ is 
$$ 0 \longrightarrow
\mathbb{Z}S  \xrightarrow{(a_{1,j})} \mathbb{Z}S^{2g} \xrightarrow{(b_{j,1})}
\mathbb{Z}S \longrightarrow \mathbb{Z} \longrightarrow 0,
 $$
where we view the $a_{1,j}$ and $b_{j,1}$ as elements of $\mathbb{Z}S$.

We take a $K(C,1)$ with one cell in each dimension such that the 
infinitely repeating exact sequence of 
left  $\mathbb{Z}C$-modules arising from the  
augmented cellular chain complex of the universal cover of the $K(C,1)$ is 
$$ 
\cdots \longrightarrow
\mathbb{Z}C \xrightarrow{me} \mathbb{Z}C \xrightarrow{c-1}  
\mathbb{Z}C \xrightarrow{me} \mathbb{Z}C \xrightarrow{c-1}
\mathbb{Z}C \longrightarrow \mathbb{Z} \longrightarrow 0.
$$

By~\cite[Theorem~3.5]{Howie05}, we get a $K(G,1)$ by melding the
one-skeleton of our $K(C,1)$ into the one-skeleton of our
$K(S,1)$ in the natural way; the attaching map of the two-cell
at the homology level is then $(a_{2,j})$.  
The exact sequence of 
left  $\mathbb{Z}G$-modules arising from the  
augmented cellular chain complex of the three-skeleton of
the universal cover of the $K(G,1)$ is 
$$ \mathbb{Z}G \xrightarrow{(0, 1-c)} 
\mathbb{Z}G^2 \xrightarrow{(a_{i,j})} \mathbb{Z}G^{2g} \xrightarrow{(b_{j,1})}
\mathbb{Z}G \longrightarrow \mathbb{Z} \longrightarrow 0.
 $$
The lemma now follows easily.
\end{proof}

We now imitate the proof of~\cite[Theorem~2]{FKS72}.

\begin{Lem}\label{L6} $G$ is of type $\VFL$ and $\chi(G) = 2-2g + \frac{1}{m}$.
\end{Lem}

\begin{proof} Let $p$ be a prime divisor of $m$. It was shown in~\cite{Baumslag62} that
$S$ is residually a finite $p$-group; see~\cite[Theorem~B]{DyerVasquez76} for an alternative proof.
Hence there exists a finite $p$-group $P = P(p)$
and a homomorphism $S \to P$ whose kernel does not contain
$\beta^{\frac{m}{p}}$, and we assume that $P$ has smallest possible order.
The centre $\Z(P)$ of $P$ is nontrivial.
By minimality of $P$, $\beta^{\frac{m}{p}}$ lies in the kernel of the composite
$S \twoheadrightarrow P \twoheadrightarrow P/\Z(P)$.
Thus $\beta^{\frac{m}{p}}$, and $\beta^{m}$, are mapped to $\Z(P)$.
By minimality of $P$, $\beta^{m}$ is mapped to $1$ in $P$. 

By considering  the direct product of such $P(p)$, one for each
prime divisor $p$ of~$m$, we find that there is a finite quotient 
of $S$ in which the image of $\beta$ has order exactly~$m$.

Hence there exists a normal subgroup $N$ of $G$ such that $N$ has finite index in $G$
and $N \cap C = \{1\}$.  It follows that $N$ acts freely on $G/C$.  The number of orbits is
$$\abs{N\backslash (G/C)} = \abs{N\backslash G/C } = \abs{(N\backslash G)/C} 
= [G:N]/m,$$
where the last equality holds  since $C$ acts freely on $N\backslash G$, on the right.

Now~\eqref{EL2} is a resolution of $\mathbb{Z}$ by
free left $\mathbb{Z}N$-modules.  Thus $N$ is of type FL, and, in particular, $N$ is
torsion-free.
It is now a simple matter to calculate $\chi(G)\,\, (= \frac{1}{[G:N]}\chi(N))$.
\end{proof}

Together  Lemma~\ref{Lg1}\eqref{Icd} and Lemma~\ref{L6} give Theorem~\ref{Tors}\eqref{Ibig1}.

By Lemma~\ref{L2}, the following is clear.
 
\begin{Cor}\label{C2} 
The sequence of left $QG$-modules
\begin{equation*}\label{Eexact}
 0 \longrightarrow
QG \oplus QGe \xrightarrow{(a_{i,j})} 
QG^{2g} \xrightarrow{(b_{j,1})}
QG \longrightarrow Q \longrightarrow 0
\end{equation*}
is exact. \hfill\qed
\end{Cor}

\begin{Lem}\label{Lcd} 
$\cd_QG=2$.
\end{Lem}

\begin{proof}  By Corollary~\ref{C2}, $\cd_QG \le 2$.  It remains to show that $\cd_QG >1$.
Let us suppose that $\cd_QG \le 1$ and derive a contradiction. 

By Notation~\ref{NH} and Lemma~\ref{L1}\eqref{I4}, 
$H$ is the (normal) subgroup of $G$ generated by the elements of finite order.  
By Dunwoody's Theorem~\cite[Theorem~IV.3.13]{DicksDunwoody89}, 
$G$ is the fundamental group of a graph of finite groups;
by~\cite[Proposition~I.7.11]{DicksDunwoody89}, $H$ is the normal 
subgroup of $G$ generated by the vertex groups.  From
the presentation of $G$ as in~\cite[Notation~I.7.1]{DicksDunwoody89}, 
it can be seen that $G/H$ is a free group.  

Since $G/H = S/\normgen{\beta}$, the abelianization of $G/H$ has $\mathbb{Z}$-rank 
$2g$ or $2g-1$. Thus the rank of the free group $G/H$ is $2g$ or $2g-1$.  
Hence $\chi(S/\normgen{\beta})$ is $1-2g$ or $2-2g$.

But $\chi(S/\normgen{\beta}) = 3-2g$ by Lemma~\ref{L6}.  This is a contradiction.
\end{proof}

Together Lemma~\ref{Lg1}\eqref{Icd}  and Lemma~\ref{Lcd} give Theorem~\ref{Tors}\eqref{Ibig2}.

By Corollary~\ref{C2} with $Q = \mathbb{C}$, the following is clear.
 
\begin{Cor}\label{C1} The  homology of
\begin{equation*}\label{EUG2}
0 \longrightarrow
\mathcal{U}(G) \oplus \mathcal{U}(G)e \xrightarrow{(a_{i,j})} 
\mathcal{U}(G)^{2g} \xrightarrow{(b_{j,1})}
\mathcal{U}(G) \longrightarrow 0
\end{equation*}
is $\Hop_*(G; \mathcal{U}(G))$. \hfill\qed
\end{Cor}

We now come to the subtle part of the argument.

\begin{Lem}\label{L3} 
 $\mathcal{U}(G) \oplus \mathcal{U}(G)e \xrightarrow{(a_{i,j})} \mathcal{U}(G)^{2g}$
is injective.
\end{Lem}

\begin{proof} Let $(u,v)$ be an arbitrary element of the kernel. Thus, 
$(u,v) \in \mathcal{U}(G) \oplus \mathcal{U}(G)e$ and
\begin{equation}\label{Eker1}
\text{ for each $j \in \{1,\ldots, 2g\}$,} \quad  u a_{1,j} +  va_{2,j} = 0 \text { in } 
\mathcal{U}(G). 
\end{equation}

Consider first the case where $u$ does not lie in $v\mathbb{C}G$.  
We shall obtain a contradiction.

We form the right  $\mathbb{C}G$-module 
$W = \mathcal{U}(G)/(v\mathbb{C}G)$, and let $w = u + v\mathbb{C}G \in W$. 
By~\eqref{Eker1}, 
\begin{equation}\label{Eker2}
\text{ for each $j \in \{1,\ldots, 2g\}$,} \quad  wa_{1,j} = 0 \text{ in } W.
\end{equation}

Let $K = \{x \in G \mid wx = w\}$.   Clearly,  $K$ is a
subgroup of $G$. 

We claim that $K  = G$; it suffices to show that
$\{x_1,\ldots, x_{2g}\} \subseteq K$.  

We will show by induction that, if $j \in \{0,1,\ldots, g\}$,
then $\{x_1,\ldots,x_{2j}\} \subseteq K$.
This is clearly true for $j = 0$.  Suppose that $j \in \{1,\ldots, g\}$
and that it is true for~$j-1$. We will show it is true for $j$.
Let $k =[x_1,x_2]\cdots[x_{2j-3},x_{2j-2}]$; then $k$ 
lies in $K$ by the induction hypothesis.  
Recall that $r_1 = [x_1,x_2]\cdots[x_{2g-1},x_{2g}]$.
By~\eqref{Eker2} and Notation~\ref{not:der},
\begin{align*} &0 = wa_{1,2j-1}= w\frac{\partial r_1}{\partial x_{2j-1}} 
= wk(1-x_{2j-1}x_{2j}x_{2j-1}^{-1})
\intertext{and}
&0 = wa_{1,2j}= w\frac{\partial r_1}{\partial x_{2j}} 
= wkx_{2j-1}(1-x_{2j}x_{2j-1}^{-1}x_{2j}^{-1}).
\end{align*}
Since $K = \{x \in G \mid w(1-x) = 0\},$ we see that
$K$ contains 
$$k(x_{2j-1}x_{2j}x_{2j-1}^{-1})k^{-1}
\text{ and } (kx_{2j-1})(x_{2j}x_{2j-1}^{-1}x_{2j}^{-1})
(kx_{2j-1})^{-1}.$$
Thus $K$ contains
$$x_{2j-1}x_{2j}x_{2j-1}^{-1} \text{ and  }
x_{2j-1}(x_{2j}x_{2j-1}^{-1}x_{2j}^{-1})x_{2j-1}^{-1},$$
and it follows easily that $K$ contains 
$x_{2j}^{-1}x_{2j-1}^{-1}$, $x_{2j-1}$ and $x_{2j}$.  This completes the
proof by induction.

Hence, $K = G$, and $w$ is fixed under the right $G$-action on  $W$.  Thus, the subset
$u + v\mathbb{C}G$ of $\mathcal{U}(G)$ is closed under the right $G$-action on $\mathcal{U}(G)$.
We denote the set $u + v\mathbb{C}G$ viewed as right $G$-set by $(u + v\mathbb{C}G)_G$.
Notice that $u + v\mathbb{C}G$ does not
contain~$0$.

By Lemma~\ref{L1}\eqref{I6}, the surjective map
$e\mathbb{C}G \to v\mathbb{C}G$, $x \mapsto vx$,
is either injective or zero.  In either event, 
$v\mathbb{C}G$ is a projective right 
$\mathbb{C}G$-module. By the left-right dual 
of~\cite[Corollary~5.6]{DicksDunwoody06}
there exists a right $G$-tree with finite edge stabilizers 
and vertex set  $(u + v\mathbb{C}G)_G$.
It follows that there exists a (left) $G$-tree $T$ with finite edge
stabilizers and vertex set 
$\null_G(u + v\mathbb{C}G)^* \subseteq \null_G(\mathcal{U}(G) -\{0\})$.

Each vertex stabilizer for $T$ is torsion,
by~\eqref{Einvertible},  
and hence embeds in $C$, by Lemma~\ref{L1}\eqref{I5}.  
By~\cite[Theorem~IV.3.13]{DicksDunwoody89},
$\cd_QG \le 1$ which contradicts Lemma~\ref{Lcd}; in essence,  $T$ is a 
one-dimensional $\uE(G)$.  Alternatively, one can use $T$ to prove that
$b_2^{(2)}(G) = 0$ and deduce that $(u,v) = (0,0)$, which is also a contradiction.

Thus $u$ lies in $v\mathbb{C}G$, and there
exists $y \in e\mathbb{C}G$ such that $u = vy$. 

We consider first the case  where $v \ne 0$. For each $j \in \{1,\ldots, 2g\}$,
$$v(ya_{1,j} + a_{2,j}) = ua_{1,j} + va_{2,j} = 0$$ by~\eqref{Eker1}, and,
by Lemma~\ref{L1}\eqref{I6}, 
$0 = ya_{1,j} + a_{2,j} = ya_{1,j} + ea_{2,j}$.  
Hence, $(y,e)$ lies in the kernel of 
$\mathbb{C}G \oplus \mathbb{C}Ge \xrightarrow{(a_{i,j})} \mathbb{C}G^{2g};$
since this map is injective by Corollary~\ref{C2}, we see $e = 0$, which is a contradiction. 

Thus $v = 0$, and hence $u = 0$.   
\end{proof}

By Lemma~\ref{L3} and Remark~\ref{Rb0} it is straightforward to obtain the following.

\begin{Lem}\label{Ldims} The $\mathcal{U}(G)$-dimensions of the kernel and the image of the map\newline
$\mathcal{U}(G) \oplus \mathcal{U}(G)e \xrightarrow{(a_{i,j})} \mathcal{U}(G)^{2g}$
are $0$ and $1 + \frac{1}{m}$, respectively.

\medskip

The $\mathcal{U}(G)$-dimensions of the image and the kernel  of  the map\newline
$\mathcal{U}(G)^{2g} \xrightarrow{(b_{j,1})} \mathcal{U}(G)$
are $1$ and $2g-1$, respectively.
 
\medskip

For $n \in \mathbb{N}$,  $b_n^{(2)}(G) = \begin{cases}
(2g-1) - (1 + \frac{1}{m}) &\text{if $n=1$,}\\
0 &\text{if $n \ne 1$.}
\end{cases}$\hfill\qed
\end{Lem}

Together Lemma~\ref{Lg1}\eqref{Ibs} and Lemma~\ref{Ldims} give Theorem~\ref{Tors}\eqref{Ibig3}.
This completes the proof of Theorem~\ref{Tors}.
\end{proof}

\section{Left-orderable groups} \label{Sorderable}

Throughout this section we will frequently make the following assumption.

\begin{Hyps}\label{Hemb}   There exist nonzero rings $Z$ and  $U$ such that  
$ZG$ is a subring of $U$ and each nonzero element of $ZG$ is invertible in $U$. 

This holds, for example, if $G$ is locally indicable, or, more generally,
left orderable, with $Z$ being any subring of $\mathbb{C}$, and $U$ being
$\mathcal{U}(G)$, by Theorem~\ref{Torder}.

Notice that  $ZG$ has no nonzero zerodivisors, and hence $G$ is torsion free. 
\hfill\qed
\end{Hyps}

\begin{Lem}\label{Lmod} 
Let $U$ be a ring,  and let  $X$ and $Y$ be sets.
 
Let $A$ and $B$  be nonzero row-finite matrices over $U$ 
in which each nonzero entry is invertible,
such that $A$ is $X \times 2$, $B$ is $2 \times Y$,  
and the  product $AB$ is the zero $X \times Y$ matrix.

Then 
$\oplus_X U \xrightarrow{A} U^2 \xrightarrow{B} \oplus_Y U$
is an exact sequence of free left $U$-modules.  

Moreover, $U^2$ has a left
$U$-basis $v_1$, $v_2$  such that  $\ker B = \im A = Uv_1$ and $B$ induces
an isomorphism $Uv_2 \simeq \im B$. 
\end{Lem}

\begin{proof} Write $A = (a_{x,i})$ and $B = (b_{i,y})$. 

There exists $x_0 \in X$ such that $(a_{x_0,1}, a_{x_0,2}) \ne (0,0)$.
We take  $v_1 = (a_{x_0,1}, a_{x_0,2})$.
Clearly $Uv_1 \subseteq \im A \subseteq \ker B$.
Without loss of generality,  there exists $y_0 \in Y$ such that 
$b_{1,y_0}$ is invertible in $U$.  We take $v_2 = (1,0)$.

Since $AB = 0$, $a_{x_0,1}b_{1,y_0} + a_{x_0,2}b_{2,y_0} = 0$.
Thus $a_{x_0,1} = - a_{x_0,2}b_{2,y_0}b_{1,y_0}^{-1}.$
Hence $a_{x_0,2}$ cannot be zero, and is therefore invertible.

Hence   $v_1$, $v_2$ is a basis of $U^2$, and 
$b_{2,y_0}b_{1,y_0}^{-1} = - a_{x_0,2}^{-1}a_{x_0,1}$.

Consider any $(a_1,a_2) \in \ker B$.
Then $a_1b_{1,y_0} + a_2 b_{2,y_0} = 0,$ and
\begin{align*}
(a_1,a_2) &= (- a_2b_{2,y_0}b_{1,y_0}^{-1},a_2) = a_2(-b_{2,y_0}b_{1,y_0}^{-1},1)
 \\&= a_2(a_{x_0,2}^{-1}a_{x_0,1},1) = a_2a_{x_0,2}^{-1}(a_{x_0,1},a_{x_0,2}) 
= a_2a_{x_0,2}^{-1}v_1 \in Uv_1,
\end{align*}
as desired.  Finally, $Uv_2 \simeq  (Uv_1 + Uv_2)/Uv_1 =  U^2/\ker B \simeq  \im B$.
\end{proof}

\begin{Rem} We see from the proof that the hypotheses that 
$A$ and $B$ are nonzero and every nonzero entry in $A$ and $B$ is invertible  
can be replaced with the hypotheses that some element of the first row of 
$B$ is invertible, and some element of the second column of $A$ is invertible.  

There are other variations, but the stated form is most convenient for our purposes.
\hfill\qed
\end{Rem}

\begin{Prop}\label{Ptwoone}  Suppose that Hypotheses~$\ref{Hemb}$ hold, and 
suppose that there exists a positive integer $n$ and a 
resolution~\eqref{Eresolution} of $Z$ by projective 
left $ZG$-modules such that  $P_n = ZG^2$.  Then  either 
the map $P_{n+1} \to P_{n}$ in~\eqref{Eresolution} is the zero map 
or  $\Hop_n(G; U) = 0$.
\end{Prop}

\begin{proof} 
 We may assume that $P_{n+1} \to P_{n}$ is nonzero.
Then we have an exact sequence 
\begin{equation}\label{EZG}
P_{n+1} \to P_n \to P_{n-1},
\end{equation}
 and we want to deduce that
\begin{equation}\label{EU}
U \otimes_{ZG} P_{n+1} \to U \otimes_{ZG}  P_n 
\to U \otimes_{ZG}  P_{n-1}
\end{equation}
remains exact.

This is clear if $ P_n \to P_{n-1}$ is the zero map.
Thus we may assume that the maps in~\eqref{EZG} are nonzero.

By adding a  suitable $ZG$-projective summand to $P_{n+1}$ 
with a zero map to $P_n$, we may assume that $P_{n+1}$ is 
$ZG$-free without affecting the images.
  Similarly, we may assume that 
 $P_{n-1}$ is $ZG$-free without affecting the kernels.  
Thus we may assume that we have specified $ZG$-bases of $P_{n+1}$,
 $P_n$ and $P_{n-1}$, and that the maps in~\eqref{EZG} 
are represented by nonzero matrices over $ZG$.

The maps in~\eqref{EU} are then represented by nonzero matrices over $U$
with all coefficients lying in $ZG$.  
Now we may apply Lemma~\ref{Lmod} to deduce that~\eqref{EU} is exact, as desired.
\end{proof}

\begin{Rem} 
 In Proposition~\ref{Ptwoone}, if we replace the hypothesis $P_n = ZG^2$ 
with the hypothesis $P_n = ZG^1$,
then it is easy to see that at least 
one of the maps $P_{n+1} \to P_n$, $P_n \to P_{n-1}$ is necessarily the zero map.
\hfill\qed
\end{Rem}

Applying Proposition~\ref{Ptwoone} with $U = \mathcal{U}(G)$,
together with Theorem~\ref{Torder}, 
we obtain the following two results.

\begin{Cor}\label{Ctwobis}  
Let $G$ be a left-orderable group, and let $Z$ be a subring of $\mathbb{C}$.
Suppose that there exists a positive integer $n$ and a resolution~\eqref{Eresolution} 
of $Z$ by projective left $ZG$-modules such that $P_n = ZG^2$.
Then either   $\cd_{Z} G \le n$ or $b_n^{(2)}(G) = 0$.
\hfill\qed
\end{Cor}

\begin{Cor}\label{Cor1}
If $G$ is a left-orderable group,
and   there exists an exact $\mathbb{C}G$-sequence of the form
\begin{equation}\label{EThompsonbis}
\cdots \xrightarrow{\partial_3} \mathbb {C}G^2
\xrightarrow{\partial_2} \mathbb{C}G^2
\xrightarrow{\partial_1} \mathbb{C}G^2
\xrightarrow{\partial_0} \mathbb{C}G
\xrightarrow{\epsilon} \mathbb{C} \longrightarrow 0
\end{equation}
in which all the $\partial_n$ are nonzero, 
then all the $b_n^{(2)}(G)$ are zero.
\end{Cor}

\begin{proof} Since $\partial_0$ is nonzero, we see that $G$ is nontrivial.
Since $G$ is torsion-free,  $b_0^{(2)}(G) = 0$.  For $n \ge 1$, 
$b_n^{(2)}(G) = 0$ by Proposition~\ref{Ptwoone}.
\end{proof}

\begin{Cor}[{L\"uck \cite[Theorem 7.10]{Lueck02}}]\label{CThompson}
All the $L^2$-Betti numbers of Thompson's group $F$ vanish.
\end{Cor}

\begin{proof} This follows from Corollary~\ref{Cor1} since 
$F$ is orderable, see~\cite{CFP96},  and 
has a resolution as in~\eqref{EThompsonbis}, 
see~\cite{BrownGeoghegan84}. 
\end{proof}

We now look at situations where we can deduce that a two-generator group is free.

\begin{Prop} \label{Pb1ne0}  Suppose that Hypotheses~$\ref{Hemb}$ hold.
The following are equivalent.

\begin{enumerate}[\normalfont(a)]
\vskip-0.6cm \null
\item
$G$ is a two-generator group, and $\Hop_1(G;U) \simeq U$.
\vskip-0.6cm \null
\item
$G$ is a two-generator group, and $\Hop_1(G;U) \ne 0$.
\vskip-0.6cm \null
\item
$G$ is free of rank two.
\end{enumerate}
\end{Prop}

\begin{proof} (a) $\Rightarrow$ (b) is obvious. 

 (b) $\Rightarrow$ (c).  Let $\{x, y\}$ be a generating set of $G$.
Then we have an exact sequence of left $ZG$-modules
$$\oplus_R ZG \longrightarrow ZG^2 
 \xrightarrow{ \begin{pmatrix}
x-1\\y-1
\end{pmatrix}} ZG  \longrightarrow Z \longrightarrow 0,$$
where $R$ is the set of relators which have a nonzero left Fox derivative
in $ZG$.  
By Proposition~\ref{Ptwoone} with $n = 1$, we see that $R$ is empty, 
and that the augmentation ideal is left $ZG$-free on $x-1$ and $y-1$.

A result of Bass-Nakayama~\cite[Proposition~1.6]{Swan69} 
then says that $G$ is freely generated by $x$ and $y$.
This can be seen geometrically, as follows. 
Let $\Gamma = \Gamma(G,\{x,y\})$ denote the Cayley graph of $G$ with respect
to the subset $\{x,y\}$.  The above exact sequence is precisely
 the augmented cellular $Z$-chain complex of~$\Gamma$.  It is then straightforward to
show that $\Gamma$ is a tree,
and that $G$ is freely generated by $x$ and $y$.

(c) $\Rightarrow$ (a) is straightforward.
\end{proof}

\begin{Cor} \label{Cb1ne0}
The following are equivalent.

\begin{enumerate}[\normalfont(a)]
\vskip-0.6cm \null
\item
$G$ is a two-generator left-orderable group and $b_1^{(2)}(G) \ne 0$.
\vskip-0.6cm \null
\item
$G$ is free of rank two. \hfill\qed
\end{enumerate}
\end{Cor}

\begin{Thm} \label{Thd1}
Suppose that Hypotheses~$\ref{Hemb}$ hold. 
If $\hd_{_{\scriptstyle Z}} G \le 1$ then every two-generator subgroup 
of $G$ is free.
\end{Thm}

\begin{proof} Since the hypotheses pass to subgroups, we may assume that
$G$ itself is generated by two elements, and it remains to show that $G$ is free.

We calculate $\Hop_*(G;U)$ in the case where $G$ is not free.

By Hypotheses~\ref{Hemb}, $G$ is torsion free.  As in Remark~\ref{Rb0}, if 
$\Hop_0(G; U) \ne 0$, 
then $G$ is free of rank zero.  Thus we may assume that $\Hop_0(G;U) = 0$.

By Proposition~\ref{Pb1ne0}, if $\Hop_1(G, U) \ne 0$, 
then  $G$ is free of rank two.
Thus we may assume that $\Hop_1(G;U) = 0$.

Since $\hd_{_{\scriptstyle Z}} G \le 1$, $\Hop_n(G;U) = 0$ 
for all $n \ge 2$.

In summary, we may assume that $\Hop_*(G; U) = 0$.

By~\cite[Theorem~4.6(b)]{Bieri81}, since $G$ is countable and $\hd_{_{\scriptstyle Z}} G \le 1$,
we have\linebreak $\cd_{_{\scriptstyle Z}} G \le 2$; in essence, the
augmentation ideal $\omega$ of $ZG$ is a countably-related flat left $ZG$-module, hence
the projective dimension of
$\null_{_{\scriptstyle ZG}}\omega$ is at most one.
Since $G$ is a two-generator group, we have a resolution of 
$Z$ by projective left 
$ZG$-modules
$$0 \longrightarrow P \longrightarrow ZG^2 
\longrightarrow ZG \longrightarrow Z \longrightarrow 0.$$
Since $\Hop_*(G;U) = 0$, 
 we have an exact sequence of projective left $U$-mod\-ules
$$0 \longrightarrow U\otimes_{ZG}P \longrightarrow  U^2 
\longrightarrow  U \longrightarrow  0.$$
This sequence splits, and we see that $\null_{_{\scriptstyle U}}(U\otimes_{\mathbb{Z}G}P)$ 
is finitely generated.

Hence $\null_{_{\scriptstyle ZG}}P$ is  finitely generated, 
by the following standard argument.  
Let $R$ be a set such that $P$ is a $ZG$-summand of $\oplus_R ZG$,
that is, $P$ is a $ZG$-submodule of $\oplus_R ZG$ and
we have a $ZG$-linear retraction of $\oplus_R ZG$ onto $P$. 
We may assume that $R$ is minimal, that is, for each $r \in R$,
 the image of $P$ under projection onto the $r$th coordinate is nonzero.
Then $U\otimes_{ZG} P$
is a $U$-submodule
of  $\oplus_R U$,
and here also $R$ is minimal.  Since  
$\null_{_{\scriptstyle U}}(U\otimes_{\mathbb{Z}G}P)$ 
is finitely generated,  
$R$ is finite, as desired.

 Now $\null_{_{\scriptstyle ZG}}Z$ has a resolution by finitely generated projective
left $ZG$-modules. By~\cite[Theorem~4.6(c)]{Bieri81}, 
$\cd_{_{\scriptstyle Z}} G \le 1$; in essence, 
$\null_{_{\scriptstyle ZG}}\omega$ is finitely related and
flat, and is therefore  projective.
Since $G$ is torsion free, $G$  is free by Stallings' Theorem; 
see Remark~II.2.3(ii) (or Corollary~IV.3.14) in~\cite{DicksDunwoody89}.
\end{proof}

\begin{Cor} \label{Chd1}
Suppose that $G$ is locally indicable, or, more generally, that $G$ is left orderable. 
If $\hd G \le 1$ then every two-generator subgroup 
of $G$ is free. \hfill\qed
\end{Cor}

We now turn from two-generator groups to  two\d1re\-la\-tor groups.

\begin{Prop} \label{Ptwotwo}  Suppose that 
$G$ is left orderable, that
$G$ has a  presentation\linebreak $\gen{X \mid R}$ with $\abs{R} = 2$, 
and that $\cd G \ge 3$.

Then 
$b_0^{(2)}(G) = 0$, $b_1^{(2)}(G) = \abs{X}-2$, and $b_2^{(2)}(G) = 0$.
\end{Prop}

\begin{proof} The given presentation of $G$ yields
an exact sequence of $\mathbb {Z}G$-modules
\[
\cdots \longrightarrow  \oplus_Y \mathbb{Z}G \xrightarrow{A} \mathbb{Z}G^2
\xrightarrow{B} \oplus_X \mathbb{Z}G
\xrightarrow{C} \mathbb{Z}G
\longrightarrow \mathbb{Z}
\longrightarrow 0.
\]
Then $H_*(G,\mathcal{U}(G))$ is the homology of the sequence
\begin{equation}\label{EUseq}
\cdots \longrightarrow  \oplus_Y \mathcal{U}(G) \xrightarrow{A} \mathcal{U}(G)^2
\xrightarrow{B} \oplus_X \mathcal{U}(G)
\xrightarrow{C} \mathcal{U}(G)
\longrightarrow 0.
\end{equation}

Since $G$ is left orderable, $G$ is torsion free.
Since $\cd G \ne 0$, $G$ is non-trivial. Hence
$G$ has an element of infinite order.
By Remark~\ref{Rb0}, $b_0^{(2)}(G) = 0$ and
 the $\mathcal{U}(G)$-dimension of  $\ker C$ in~\eqref{EUseq} is $\abs{X}-1$.

Since $G$ is left orderable, all 
nonzero elements of $\mathbb {C}G$ are invertible in
$\mathcal{U}(G)$ by Theorem \ref{Torder}.  Since $\cd G \ge 3$, 
 $b_2^{(2)}(G) = 0$ by Corollary~\ref{Ctwobis}.   
Moreover, by Lemma~\ref{Lmod}, the $\mathcal{U}(G)$-dimension of 
$\im B$  in~\eqref{EUseq} is one. 

Finally, $b_1^{(2)}$ is the difference between 
the $\mathcal{U}(G)$-dimensions of  $\ker C$ and  $\im B$ in~\eqref{EUseq}, 
that is, $\abs{X}-2$.  Of course, the hypotheses clearly imply that $\abs{X} \ge 2$.
\end{proof}

Suppose that $G$ is a left-orderable two-relator group.
We know the first three $L^2$-Betti numbers of $G$ if 
$\cd G \ge 3$ by Proposition~\ref{Ptwotwo}.
If $\cd G \le 1$, then $G$ is free, and
again one knows the $L^2$-Betti numbers.
There remains the case where $\cd G = 2$;
here all we know are the $L^2$-Betti numbers of
torsion-free \ors groups; these groups are left-orderable
by~\cite[Theorem~2.2]{Hempel90} and they are clearly two-relator groups.

\bigskip

\noindent{\textbf{\Large{Acknowledgments}}}

\medskip
\footnotesize

The research of the first-named author  was
funded by the DGI (Spain) through Project  BFM2003-06613.

\bibliographystyle{amsplain}

\begin{thebibliography}{21}
\bibitem{Baumslag62}
Gilbert Baumslag,
\newblock{\em   On generalized free products},
\newblock Math.\ Z.\ \textbf{78}(1962), 423--438.
\vskip-0.4cm \null
\bibitem{Bieri81}
 Robert Bieri,
\newblock {\em Homological dimension of discrete groups, second edition},
\newblock Queen Mary College Mathematical Notes,  London,  1981.
\vskip-0.4cm \null
\bibitem{Brodskii80}
S.\ D.\ Brodski\u{\i},
\newblock{\em  Equations over groups and groups with one defining relation},
\newblock Siberian Math J.\ \textbf{25}(1984), 235--251.
\vskip-0.4cm \null
\bibitem{BrownGeoghegan84}
Kenneth S.\ Brown and Ross Geoghegan,
\newblock {\em An infinite-dimensional torsion-free {${\rm FP}\sb{\infty }$} group},
\newblock Invent.\ Math.\ \textbf{77}(1984), 367--381.
\vskip-0.4cm \null
\bibitem{BurnsHale72}
R.\ G.\ Burns and V.\ W.\ D.\ Hale,
\newblock {\em A note on group rings of certain torsion-free groups},
\newblock Canad.\ Math.\ Bull.\ \textbf{15}(1972), 441--445.
\vskip-0.4cm \null
\bibitem{CFP96}
J.\ W.\ Cannon, W.\ J.\ Floyd and W.\ R.\ Parry,
\newblock {\em Introductory notes on {R}ichard {T}hompson's groups},
\newblock Enseign.\ Math.\ \textbf{42}(1996), 215--256.
\vskip-0.4cm \null
\bibitem{Chiswell03}
I.\ M.\ Chiswell,
\newblock {\em Euler characteristics of discrete groups},
\newblock pp. 106--254 in: Groups: topological, combinatorial
  and arithmetic aspects (ed. T.~W. M{\"u}ller),   
London Math.\ Soc.  Lecture Note Ser.\ \textbf{311}, 
CUP, Cambridge, 2004.
\vskip-0.4cm \null
\bibitem{DicksDunwoody89}
Warren Dicks and M.\ J.\ Dunwoody,
\newblock {\em Groups acting on graphs}, 
\newblock Cambridge Stud. Adv. Math.  \textbf{17}, CUP,  Cambridge, 1989.
\newline  Errata at: http://mat.uab.cat/$\scriptstyle\sim$dicks/DDerr.html.
\vskip-0.4cm \null
\bibitem{DicksDunwoody06}
Warren Dicks and M.\ J.\ Dunwoody,
\newblock {\em Retracts of vertex sets of trees and the almost stability theorem}, 
\newblock  math.GR/0510151 version 3
\vskip-0.4cm \null
\bibitem{DyerVasquez76}
Eldon Dyer and A.\ T.\ Vasquez, 
\newblock{\em Some properties of two-dimensional Poincar\'e duality groups},
\newblock pp. 45--54 in: 
Algebra, topology, and category theory:
a collection of papers in honor of Samuel Eilenberg 
(eds. Alex Heller,  Myles Tierney and Samuel Eilenberg),
Academic Press, New York, 1976.
\vskip-0.4cm \null
\bibitem{FKS72}
J.\ Fischer, A.\ Karrass and D.\ Solitar,
\newblock {\em On \ore  groups having elements of finite order},
\newblock Proc.\ Amer.\ Math.\ Soc.\ \textbf{33}(1972), 297--301.
\vskip-0.4cm \null
\bibitem{Hempel90}
John Hempel,
\newblock {\em \Ore surface groups},
\newblock Math.\ Proc.\ Cambridge Philos.\ Soc. \textbf{108}(1990), 467--474.
\vskip-0.4cm \null
\bibitem{Howie82}
James Howie,
\newblock {\em On locally indicable groups},
\newblock Math.\ Z.\ \textbf{180}(1982),  445--461.
\vskip-0.4cm \null
\bibitem{Howie00}
James Howie,
\newblock {\em A short proof of a theorem of {B}rodski\u{\i}},
\newblock Publ.\ Mat.\ \textbf{44}(2000), 641--647.
\vskip-0.4cm \null
\bibitem{Howie05}
James Howie,
\newblock {\em Some results on \ore  surface groups},
\newblock Bol. Soc. Mat. Mexicana (3), \textbf{10} (Special
Issue) (2004), 255--262.  Erratum ibid, 545--546.
\vskip-0.4cm \null
\bibitem{KLL04}
Peter Kropholler, Peter Linnell and Wolfgang L{\"u}ck,
\newblock {\em Groups of small homological dimension and the {A}tiyah conjecture},
\newblock math.GR/0401312.
\vskip-0.4cm \null
\bibitem{Linnell92}
Peter A.\ Linnell,
\newblock {\em Zero divisors and {$L\sp 2(G)$}},
\newblock C.\ R.\ Acad.\ Sci.\ Paris S\'er.\ I Math.\ \textbf{315}(1992), 49--53.
\vskip-0.4cm \null
\bibitem{Linnell93}
Peter A.\ Linnell,
\newblock {\em Division rings and group von {N}eumann algebras},
\newblock Forum Math.\ \textbf{5}(1993), 561--576.
\vskip-0.4cm \null
\bibitem{Lueck02}
Wolfgang L{\"u}ck,
\newblock {\em {$L\sp 2$}-invariants: theory and applications to geometry and
  {$K$}-theory}, 
\newblock{Ergeb. Math. Grenzgeb.(3) \textbf{44}, Springer-Verlag, Berlin, 2002.}
\vskip-0.4cm \null
\bibitem{Richards63}
Ian Richards,
\newblock{\em On the classification of noncompact surfaces},
\newblock Trans.\ Amer.\ Math.\ Soc.\  \textbf{106}(1963), 259--269.
\vskip-0.4cm \null
\bibitem{Swan69}
Richard G.\ Swan,
\newblock {\em Groups of cohomological dimension one}, 
\newblock J.\ Algebra \textbf{12}(1969), 585--601.
\end{thebibliography}

\textsc{Departament de  Matem\`atiques,
Universitat Aut\`onoma de Barcelona,
E-08193 Bellaterra (Barcelona), Spain}

\emph{E-mail address}{:\;\;}\texttt{dicks@mat.uab.cat}

\emph{URL}{:\;\;}\texttt{http://mat.uab.cat/$\scriptstyle\sim$dicks/}

\medskip

\textsc{Department of Mathematics, 
Virginia Tech, 
Blacksburg, VA 24061-0123, USA}

\emph{E-mail address}{:\;\;}\texttt{linnell@math.vt.edu}

\emph{URL}{:\;\;}\texttt{http://www.math.vt.edu/people/linnell/}
\end{document}